\input amstex
\input amsppt.sty
\magnification=\magstep1
\vsize=22.2truecm
\baselineskip=16truept
\nologo
\pageno=1
\TagsOnRight
\def\Z{\Bbb Z}
\def\N{\Bbb N}

\def\l{\left}
\def\r{\right}
\def\bg{\bigg}
\def\({\bg(}
\def\[{\bg[}
\def\){\bg)}
\def\]{\bg]}
\def\t{\text}
\def\f{\frac}

\def\bi{\binom}
\def\eq{\equiv}

\def\ls{\leqslant}
\def\gs{\geqslant}

\def\ve{\varepsilon}
\def\da{\delta}

\def\bi{\binom}

\def\Proof{\noindent{\it Proof}}
\def\Remark{\medskip\noindent{\it  Remark}}
\def\Ack{\medskip\noindent {\bf Acknowledgments}}
\def\pmod #1{\ (\roman{mod}\ #1)}
\def\mo{\roman{mod}}

\topmatter \hbox{Discrete Math. 306(2006), no.\,16, 1921--1940.}
\bigskip
\title {A combinatorial identity with application to Catalan numbers}\endtitle
\author {Hao Pan and Zhi-Wei Sun}\endauthor
\affil Department of Mathematics, Nanjing University
\\Nanjing 210093, People's Republic of China
\\ {\tt haopan79\@yahoo.com.cn}\ \ \quad\tt{zwsun\@nju.edu.cn}
\endaffil
\abstract By a very simple argument, we prove that if $l,m,n\in\{0,1,2,\ldots\}$ then
$$\sum_{k=0}^l(-1)^{m-k}\bi{l}{k}\bi{m-k}{n}\bi{2k}{k-2l+m}
=\sum_{k=0}^l\bi lk\bi{2k}n\bi{n-l}{m+n-3k-l}.$$
On the basis of this identity, for $d,r\in\{0,1,2,\ldots\}$
we construct explicit $F(d,r)$ and $G(d,r)$ such that
for any prime $p>\max\{d,r\}$ we have
$$\sum_{k=1}^{p-1}k^{r}C_{k+d}\eq\cases F(d,r)\ (\mo\ p)&\t{if}\ p\eq1\ (\mo\ 3),
\\G(d,r)\ (\mo\ p)&\t{if}\ p\eq 2\ (\mo\ 3),\endcases$$
where $C_n$ denotes the Catalan number $\f1{n+1}\bi{2n}n$. For example,
when $p\gs 5$ is a prime, we have
$$\sum_{k=1}^{p-1}k^2C_k\eq\cases-2/3\ (\mo\ p)&\t{if}\ p\eq1\ (\mo\ 3),
\\-1/3\ (\mo\ p)&\t{if}\ p\eq2\ (\mo\ 3);\endcases$$
and
$$\sum_{0<k<p-4}\f{C_{k+4}}k
\eq\cases 503/30\ (\mo\ p)&\t{if}\ p\eq1\ (\mo\ 3),\\-100/3\ (\mo\ p)&\t{if}\ p\eq2\ (\mo\ 3).
\endcases$$
This paper also contains some new recurrence relations for Catalan numbers.
\endabstract
\keywords {Binomial coefficient; Combinatorial identity; Catalan number}
\endkeywords
\thanks  2000 {\it Mathematics Subject Classification}.
Primary 11B65; Secondary 05A10, 05A19, 11A07, 11B37.
\newline\indent The second author is responsible for communications,
and partially supported by the National Science Fund for
Distinguished Young Scholars (no. 10425103) and a Key Program of
NSF in China.
\endthanks
\endtopmatter

\document

\heading{1. Introduction}\endheading

As usual, for $k\in\Z$ we define the binomial coefficient $\bi xk$ as follows:
$$\bi xk=\cases \f1{k!}\prod_{j=0}^{k-1}(x-j)&\t{if}\ k>0,
\\1&\t{if}\ k=0,
\\0&\t{if}\ k<0.\endcases$$
There are many combinatorial identities involving binomial coefficients.
(See, e.g., [GJ], [GKP] and [PWZ].) A nice identity of Dixon (cf. [PWZ, p.\,43]) states that
$$\sum_{k\in\Z}(-1)^k\bi{a+b}{a+k}\bi{b+c}{b+k}\bi{c+a}{c+k}=\f{(a+b+c)!}{a!b!c!}$$
for any $a,b,c\in\N=\{0,1,2,\ldots\}$.

During the second author's visit (January--March, 2005)
to the Institute of Camille Jordan at Univ. Lyon-I,
Dr. Victor J. W. Guo told Sun that he had made the following ``conjecture":
{\it Given $l,m\in\N$ one has
$$\sum_{k=0}^l(-1)^{m-k}\bi{l}{k}\bi{m-k}{l}\bi{2k}{k-2l+m}=\cases
\binom{2m/3}{m/3}\binom{m/3}{l-m/3}&\text{if}\ \ 3\mid m,\\
0&\text{otherwise};
\endcases$$
in other words,}
$$\sum_{k=0}^l(-1)^{m-k}\bi{l}{k}\bi{m-k}{l}\bi{2k}{k-2l+m}=[3\mid m]\bi l{\lceil m/3\rceil}
\bi{2\lceil m/3\rceil}l,\tag1.0$$
where $\lceil\cdot\rceil$ is the ceiling function, and for
an assertion $A$ we adopt the notation
$$[A]=\cases1&\t{if}\ A\ \t{holds},\\0&\t{otherwise}.\endcases$$

The above conjecture is similar to Dixon's identity in some sense; of course
it can be proved with the aid of computer
via the WZ method or Zeilberger's algorithm (cf. [PWZ]).
After we showed (1.0)
in a preliminary version of this paper by
Lagrange's inversion formula (cf. [GJ, p.\,17]),
Prof. C. Krattenthaler at Univ. Lyon-I kindly told us that (1.0)
can also be proved by letting $a=m-3l$, $b=1/2-l$ and $x\to 1$ in Bailey's hypergeometric
series identity (cf. [B] or Ex. 38(a) of [AAR, p.\,185])
$$\align
&{}_3F_2\(\matrix a,\, 2b-a-1,\, a-2b+2\\b,\, a-b+3/2\endmatrix\,;\,\f x4\)\\
=&\f1{(1-x)^{a}}\ {}_3F_2\(\matrix a/3,\, (a+1)/3,\, (a+2)/3
\\b,\, a-b+3/2\endmatrix\,;\,-\f{27x}{4(1-x)^3}\).
\endalign$$

In this paper, by a simple argument
we show the following combinatorial identity
the special case $n=l$ of which yields (1.0).

\proclaim{Theorem 1.1} Provided that $l,m,n\in \N$, we have
$$\sum_{k=0}^l(-1)^{m-k}\bi{l}{k}\bi{m-k}{n}\bi{2k}{k-2l+m}
=\sum_{k=0}^l\bi lk\bi{2k}n\bi{n-l}{m+n-3k-l}.\tag1.1$$
\endproclaim

\Remark\ 1.1. (a) The preceding hypergeometric
series identity of Bailey does not imply (1.1) which involves
three parameters $l,\,m$ and $n$. However, Prof. C. Krattenthaler
informed us that (1.1) can also be deduced by putting
$a=m/3-l$, $b=d=1-2l+m$ and $e=1-l+m-n$ in the complicated hypergeometric
identity (3.26) of [KR] (which was obtained on the basis of Bailey's identity).
Nevertheless, (1.1) has not been pointed out explicitly before, and
our proof of (1.1) is very elementary and particularly
simple.

(b) The identity (1.1) might have a combinatorial interpretation
related to Callan's idea (cf. [C]) in his combinatorial proof
of a curious identity due to Sun.

\proclaim{Corollary 1.1} Let $l$ and $m$ be nonnegative integers. Then
$$\sum_{k=0}^l(-1)^{m-k}\bi lk\bi{m-k}{l+1}\bi{2k}{k-2l+m}
=(1-[3\mid m-1])\bi{l}{\lceil m/3\rceil}
\bi{2\lceil m/3\rceil}{l+1}\tag1.2$$
and $$\sum_{k=0}^l(-1)^{m-k}\bi lk\bi{m-k}{l+2}\bi{2k}{k-2l+m}
=(1+[3\mid m+1])\bi{l}{\lceil m/3\rceil}\bi{2\lceil m/3\rceil}{l+2}.\tag1.3$$
\endproclaim
\Proof. Putting $n=l+j$ in (1.1) with $j\in\{1,2\}$, we get that
$$\sum_{k=0}^l(-1)^{m-k}\bi lk\bi{m-k}{l+j}\bi{2k}{k-2l+m}
=\sum_{k=0}^l\bi lk\bi{2k}{l+j}\bi j{j-(3k-m)}.$$
If $0\ls 3k-m\ls j$, then $m/3\ls k\ls (m+2)/3$ and hence $k=\lceil m/3\rceil$.
Note that
$$\bi 1{1-(3\lceil m/3\rceil-m)}=1-[3\mid m-1]\ \ \t{and}\ \
\bi 2{2-(3\lceil m/3\rceil-m)}=1+[3\mid m+1].$$
So we have (1.2) and (1.3). \qed
\smallskip

From (1.0), (1.2) and (1.3) we can deduce the following result.

\proclaim{Theorem 1.2} Let $p$ be a prime and $d\in\{0,\ldots,p\}$. Then
$$\sum_{k=0}^{p-1}\bi{2k}{k+d}\eq\l(\f{p-d}3\r)\pmod{p},\tag1.4$$
where the Legendre symbol $(\f a3)$ coincides with the unique integer in $\{0,\pm1\}$
satisfying $a\eq(\f a3)\pmod{3}$.
Also,
$$\sum_{k=1}^{p-1}k\bi{2k}{k+d}
\eq\l([3\mid p-d]-\f13\r)\l(2\l(\f{p-d}3\r)-d\r)-[p=3]\pmod{p},\tag1.5$$
and
$$\sum_{k=1}^{p-1}\f{\bi{2k}{k+d}}k\equiv
\cases d^{-1}(-1+2(-1)^d+3[3\mid p-d])\pmod{p}&\t{if}\ d\not=0,
\\-[p=3]\pmod{p}&\t{if}\ d=0.\endcases\tag1.6$$
\endproclaim

The well-known Catalan numbers given by
$$C_n=\f1{n+1}\bi{2n}n=\bi{2n}n-\bi{2n}{n-1}\qquad (n=0,1,2,\ldots)$$
play important roles in combinatorics. For $n\in\N$ and $j=0,1,\ldots,n+1$,
we define
$$C_{n,j}=2\bi{2n}{n-j}-\bi{2n}{n-1-j}-\bi{2n}{n+1-j}$$
and view $C_{n,j}/2$ as a generalized Catalan number;
it is clear that $C_{n,0}/2=C_n$.
From (1.6) we can deduce the following result.

\proclaim{Corollary 1.2} Let $p$ be a prime. Then, for any $d=0,\ldots,p-1$ we have
$$\sum_{k=1}^{p-1}\f{C_{k+d}}k\eq
-[p=3]C_d+\sum_{j=1}^{d+1}\l(-1+2(-1)^j+3[3\mid p-j]\r)\f{C_{d,j}}j
\pmod{p}.\tag1.7$$
Consequently, if $p\gs5$ then
$$\align
\\&\sum_{k=1}^{p-1}\f{C_k}k\eq\f32\l(1-\l(\f p3\r)\r)\ (\mo\ p),\tag1.8
\\&\sum_{k=1}^{p-2}\f{C_{k+1}}k\eq\f34\l(1+\l(\f p3\r)\r)\ (\mo\ p),\tag1.9
\\&\sum_{k=1}^{p-3}\f{C_{k+2}}k\eq3\l(\f p3\r)\ (\mo\ p),\tag1.10
\\&\sum_{k=1}^{p-4}\f{C_{k+3}}k\eq\f{207(\f p3)-47}{24}\ (\mo\ p),\tag1.11
\\&\sum_{k=1}^{p-5}\f{C_{k+4}}k\eq\f{1503(\f p3)-497}{60}\ (\mo\ p).\tag1.12
\endalign$$
\endproclaim
\Proof. Let $d\in\{0,\ldots,p-1\}$ and $k\in\N$.
With the help of the Chu-Vandermonde identity (cf. [GKP, (5.27)]),
$$\align C_{k+d}=&\bi{2k+2d}{k+d}-\bi{2k+2d}{k+d-1}
\\=&\sum_{j=-d}^d\bi{2d}{d-j}\bi{2k}{k+j}-\sum_{j=-d}^d\bi{2d}{d-j}\bi{2k}{k+j-1}
\\=&\sum_{j=-d}^d\bi{2d}{d-j}\bi{2k}{k+j}-\sum_{i=-d-1}^{d+1}\bi{2d}{d-1-i}\bi{2k}{k+i}
\\=&\sum_{0<j\ls d}\(\bi{2d}{d-j}\bi{2k}{k+j}+\bi{2d}{d+j}\bi{2k}{k-j}\)
\\&-\sum_{0<j\ls d+1}\(\bi{2d}{d-1-j}\bi{2k}{k+j}+\bi{2d}{d-1+j}\bi{2k}{k-j}\)
\\&+\bi{2d}d\bi{2k}k-\bi{2d}{d-1}\bi{2k}k.
\endalign$$
Thus
$$\align&C_{k+d}+\(\bi{2d}{d-1}-\bi{2d}d\)\bi{2k}k
\\=&\sum_{0<j\ls d}2\bi{2d}{d-j}\bi{2k}{k+j}
-\sum_{j=1}^{d+1}\(\bi{2d}{d-1-j}+\bi{2d}{d+1-j}\)\bi{2k}{k+j}
\endalign$$
and hence
$$C_{k+d}=C_d\bi{2k}k+\sum_{j=1}^{d+1}C_{d,j}\bi{2k}{k+j}.\tag1.13$$

In view of (1.13),
$$\sum_{k=1}^{p-1}\f{C_{k+d}}k=C_d\sum_{k=1}^{p-1}\f{\bi{2k}k}k
+\sum_{j=1}^{d+1}C_{d,j}\sum_{k=1}^{p-1}\f{\bi{2k}{k+j}}k.$$
Combining this with (1.6), we immediately get (1.7).

Observe that
$$C_{p+k}\eq2C_k\ (\mo\ p)\qquad \ \t{for every}\ k=0,\ldots,p-2;\tag1.14$$
in fact,
$$C_{p+k}=\f{\bi{2p+2k}{p+k}}{p+k+1}
=\f{\bi{2p}p}{p+k+1}\cdot\f{\prod_{0<j\ls 2k}(2p+j)}{\prod_{0<i\ls k}(p+i)^2}
\eq\f2{k+1}\bi{2k}k\ (\mo\ p)$$
since $$\f12\bi{2p}p=\bi{2p-1}{p-1}=\prod_{j=1}^{p-1}\f{2p-j}j\eq(-1)^{p-1}\eq1\ (\mo\ p).$$
Thus
$$\sum_{p-d\ls k<p}\f{C_{k+d}}k=\sum_{0<k\ls d}\f{C_{p-k+d}}{p-k}
\eq-2\sum_{0<k\ls d}\f{C_{d-k}}k\ (\mo\ p).$$
(Note that if $0<k\ls d$ then $0\ls d-k<d\ls p-1$.)

Now assume that $p\gs 5$. Clearly
$$\align\sum_{j=1}^{d+1}[3\mid{p-j}]\f{C_{d,j}}j
=&\f{1+(\f p3)}2\sum\Sb 1\ls j\ls d+1\\j\eq1\,(\mo\ 3)\endSb\f{C_{d,j}}j
+\f{1-(\f p3)}2\sum\Sb 1\ls j\ls d+1\\j\eq2\,(\mo\ 3)\endSb\f{C_{d,j}}j
\\=&\f12\sum^{d+1}\Sb j=1\\3\nmid j\endSb\f{C_{d,j}}j
+\l(\f p3\r)\f12\sum_{j=1}^{d+1}\l(\f j3\r)\f{C_{d,j}}j.
\endalign$$
Therefore, by applying the above and (1.7) we obtain that
$$\align&\sum_{0<k<p-d}\f{C_{k+d}}k
-2\sum_{0<k\ls d}\f{C_{d-k}}k
\\\eq&\sum_{k=0}^{p-1}\f{C_{k+d}}k
\eq\sum_{j=1}^{d+1}\l(2(-1)^j-1\r)\f{C_{d,j}}j
+\f32\sum^{d+1}\Sb j=1\\3\nmid j\endSb\f{C_{d,j}}j
+\l(\f p3\r)\f32\sum_{j=1}^{d+1}\l(\f j3\r)\f{C_{d,j}}j\ (\mo\ p).
\endalign$$
When $d=0,1,2,3,4$, this yields (1.8)--(1.12)
after some trivial computations. \qed
\smallskip

As usual we let $\lfloor\cdot\rfloor$ be the greatest integer function.
On the basis of Theorem 1.1, we also establish the following general theorem
concerning Catalan numbers.

\proclaim{Theorem 1.3} Let $p$ be a prime and $d,r\in\{0,\ldots,p-1\}$.
Then
$$(-1)^r\sum_{k=0}^{p-1}\bi{k+r}rC_{k+d}\eq\sum_{0\ls k<d}\bi{d-1-k}rC_k
+\sum_{i=0}^{r}(-1)^{i}\bi d{r-i}f_i(\ve_i)\ (\mo\ p),\tag1.15$$
where $\ve_i=(\f{p-i-1}3)$ and
$$\align f_i(\ve_i)=&\sum_{k=0}^{\lfloor (i+1-\ve_i)/3\rfloor}
(-1)^{k+\ve_i}\bi{i+2}{3k+1+\ve_i}\bi{k+(i-2+\ve_i)/3}i
\\&+[\ve_i=0\ \&\ 3\mid i+1]+[i=0](3[\ve_i=-1]-1).
\endalign$$
\endproclaim

\Remark\ 1.2. (a) (1.15) in the case $d=0$
yields the congruence
$$\aligned &\sum_{k=0}^{p-1}\bi{k+r}rC_k
\\\eq&\cases\sum_{k=0}^{\lfloor(r+2)/3\rfloor}
(-1)^{k-1}\bi{r+2}{3k}\bi{k+(r-3)/3}r\ (\mo\ p)&\t{if}\ p-r\eq0\ (\mo\ 3),
\\\sum_{k=0}^{\lfloor(r+1)/3\rfloor}
(-1)^k\bi{r+2}{3k+1}\bi{k+(r-2)/3}r+[p=3]\ (\mo\ p)&\t{if}\ p-r\eq1\ (\mo\ 3),
\\\sum_{k=0}^{\lfloor r/3\rfloor}
(-1)^{k-1}\bi{r+2}{3k+2}\bi{k+(r-1)/3}r\ (\mo\ p)&\t{if}\ p-r\eq2\ (\mo\ 3).
\endcases
\endaligned$$

(b) Let $p$ be a prime and $d\in\{0,\ldots,p-1\}$.
For each $r=0,\ldots,p-1$, clearly
$$\align(-1)^r\sum_{k=1}^pk^rC_{k+d-1}=&\sum_{k=0}^{p-1}(-k-1)^rC_{k+d}
=\sum_{k=0}^{p-1}C_{k+d}\sum_{s=0}^rs!S(r,s)\bi{-k-1}s
\\=&\sum_{s=0}^r(-1)^ss!S(r,s)\sum_{k=0}^{p-1}\bi{k+s}sC_{k+d},
\endalign$$ where
$$S(r,s)=\f1{s!}\sum_{t=0}^s(-1)^{s-t}\bi stt^r\ \ \ (0\ls s\ls r)$$
are Stirling numbers of the second kind (cf. [GKP]).
This, together with Theorem 1.3, shows that if $P(x)$ is a polynomial
of degree at most $p-1$ with $p$-adic integer
coefficients then
$$\sum_{k=0}^{p-1}P(k)C_{k+d}\eq \psi_P(d;\,p\ \mo\ 3)\ \ (\mo\ p)$$
for a suitable function $\psi_P$ which can be constructed explicitly.
This is general enough, because any integer $r$ can be written in the form
$(p-1)q+r_0$ with $q\in\Z$ and $r_0\in\{0,\ldots,p-2\}$,
and by Fermat's little theorem
we have $k^{r}\eq k^{r_0}\ (\mo\ p)$ for all $k=1,\ldots,p-1$.
\smallskip

\proclaim{Corollary 1.3} Let $p$ be a prime, and let $d\in\{0,1,\ldots,p-1\}$. Then we have
$$\gather\sum_{k=0}^{p-1}C_{k+d}\eq\f{3(\f p3)-1}2+\sum_{0\ls k<d}C_k\ (\mo\ p),\tag1.16
\\\sum_{k=0}^{p-1}kC_{k+d}\eq\f{d+1}2\l(1-\l(\f p3\r)\r)-\l(\f p3\r)d
-\sum_{k=0}^dkC_{d-k}\ (\mo\ p),\tag1.17
\\\sum_{k=0}^{p-1}k^2C_{k+d}
\eq\f{9d^2+6d-1}6\l(\f p3\r)-\f{(d+1)^2}2-[p=3]+\sum_{0<k\ls d}k^2C_{d-k}\ (\mo\ p).\tag 1.18
\endgather$$
\endproclaim
\Proof. For $i=0,1,2$ let $\ve_i$ and $f_i(\ve_i)$ be as in Theorem 1.3.
It is easy to verify that
$$f_0(\ve_0)=(-1)^{\ve_0}\bi 2{1+\ve_0}+3[\ve_0=-1]-1\eq\f{3(\f p3)-1}2\ (\mo\ p)$$
and $f_1(\ve_1)=(\f p3)$ and $f_2(\ve_2)=\f 23(\f p3)+[p=3]$.
Thus (1.15) in the case $r=0$ is actually equivalent to (1.16).
Putting $r=1$ in (1.15) we get that
$$\align-\sum_{k=0}^{p-1}(k+1)C_{k+d}\eq&\sum_{0\ls k<d}(d-1-k)C_k+df_0(\ve_0)-f_1(\ve_1)
\\\eq&-\sum_{0\ls k<d}C_k+\sum_{k=0}^d(d-k)C_k+d\f{3(\f p3)-1}2-\l(\f p3\r)\ (\mo\ p).
\endalign$$
This, together with (1.16), yields (1.17).
By (1.15) in the case $r=2$,
$$\align\sum_{k=0}^{p-1}\f{(k+1)(k+2)}2C_{k+d}
\eq&\sum_{0\ls k<d}\f{(d-k-1)(d-k-2)}2C_k
\\&+\f{d(d-1)}2f_0(\ve_0)-df_1(\ve_1)+f_2(\ve_2)\ (\mo\ p)
\endalign$$
and hence
$$\align\sum_{k=0}^{p-1}(k^2+3k+2)C_{k+d}
\eq&\sum_{0\ls k<d}((d-k)^2-3(d-k)+2)C_k
+d(d-1)\f{3(\f p3)-1}2
\\&-2d\l(\f p3\r)+\f 43\l(\f p3\r)+2[p=3]\ \ (\mo\ p).
\endalign$$
Combining this with (1.16) and (1.17) we immediately get (1.18). \qed

\smallskip

The Catalan numbers can also be defined
by $C_0=1$ and the recursion $C_{n+1}=\sum_{k=0}^nC_kC_{n-k}\ (n=0,1,2,\ldots)$.
Below we provide some new recursions for Catalan numbers by using our previous congruences.

\proclaim{Theorem 1.4} Let $d\in\N$ and $\da\in\{0,1\}$. Then we have
$$\aligned C_d=&(1-2\da)\sum_{0\ls k<d}C_k
+(-1)^{\delta}\sum_{i=0}^d\l(\f{i-\da}3\r)C_{d,i+1}+1+\da
\\=&\f12\sum_{j=1}^{d+1}(1-3[3\mid j])C_{d,j}+\f 32.\endaligned\tag1.19$$
Also,
$$\aligned\sum_{k=0}^dkC_{d-k}
=&\sum_{j=1}^{d+1}\l(\f{j-1}3\r)\(2\bi{2d}{d-j}-(d+1)C_{d,j}\)-d
\\=&\sum_{j=1}^{d+1}\l(\f{j+1}3\r)\(2\bi{2d}{d-j}-(d+1)C_{d,j}\)+2d+1
\endaligned\tag1.20$$
and
$$\aligned&\sum_{0<k\ls d}\l(k-\f 23\r)C_{d-k}+\f13\sum_{j=1}^{d+1}jC_{d,j}
\\=&\sum\Sb 1\ls j\ls d+1\\j\eq1\,(\mo\ 3)\endSb jC_{d,j}-d+\f 23
=\sum\Sb 1\ls j\ls d+1\\j\eq2\,(\mo\ 3)\endSb jC_{d,j}+2d-\f 13.
\endaligned\tag1.21$$
\endproclaim

\Remark\ 1.3. A referee of this paper noted that some identities in Theorem 1.4,
such as (1.19), can also be established by generating function manipulations
and the observation
$$C_{d,j}=C_{d-j}^{(2j)}-C_{d-(j-1)}^{(2j-2)}\ \ \ \ (j=1,\ldots,d+1),$$
where $C_n^{(k)}=\bi{2n+k}n-\bi{2n+k}{n-1}$ for $k,n\in\Z$.

\medskip

In Sections 2-5 we are going to show Theorems 1.1-1.4 respectively.

\heading{2. Proof of Theorem 1.1}\endheading

Let $R$ be a commutative ring with identity.
For a formal power series
$f(t)\in R[\![t]\!]$ and a nonnegative integer $n$,
by $[t^n]f(t)$ we mean the coefficient of
$t^n$ in $f(t)$.

\medskip
\noindent{\it Proof of Theorem 1.1}.
We fix $l,n\in\N$.
By the Chu-Vandermonde identity, we have
$$\align
&\sum_{m=0}^{\infty}s^m\sum_{k=0}^l(-1)^{m-k}\bi{l}{k}\bi{m-k}{n}\bi{2k}{k-2l+m}\\
=&\sum_{m=0}^{\infty}s^m\sum_{k=0}^l(-1)^{m-k}\bi{l}{k}\sum_{j=0}^n\bi{2(l-k)}{n-j}
\bi{k-2l+m}{j}\bi{2k}{k-2l+m}\\
=&\sum_{k=0}^l(-1)^k\bi{l}{k}\sum_{j=0}^n\bi{2(l-k)}{n-j}\bi{2k}{j}
\sum_{m=0}^{\infty}\bi{2k-j}{k-2l+m-j}(-s)^m\\
=&\sum_{k=0}^l(-1)^k\bi{l}{k}\sum_{j=0}^n\bi{2(l-k)}{n-j}
\bi{2k}{j}(-s)^{2l-k+j}(1-s)^{2k-j}
\endalign$$
and hence
$$\align
&\sum_{m=0}^{\infty}s^m\sum_{k=0}^l(-1)^{m-k}\bi{l}{k}\bi{m-k}{n}\bi{2k}{k-2l+m}\\
=&\sum_{k=0}^l(-1)^k\bi{l}{k}(-s)^{2l-k}
[t^n](1+t)^{2(l-k)}((1-s)-st)^{2k}\\
=&[t^n]s^{l}\sum_{k=0}^l\bi{l}{k}\l(s(1+t)^2\r)^{l-k}
(1-s(1+t))^{2k}
\\=&[t^n]s^{l}\l(s(1+t)^2+(1-s(1+t))^2\r)^l
\\=&[t^n]\l((s+st)^2+s(1-s-st)^2\r)^l.
\endalign$$

Clearly
$$\align&[t^n](1+s)^l\l((s+t)^2+s(1-s-t)^2\r)^l
\\=&[t^n]\l((1+s)((1+s)t^2+2s^2t+s^2+s(1-s)^2)\r)^l
\\=&[t^n]\l(((1+s)t+s^2)^2+s\r)^l
\\=&\sum_{k=0}^l\bi lk[t^n]\l((1+s)t+s^2\r)^{2k}s^{l-k}
\\=&\sum_{k=0}^l\bi lk\bi{2k}n(1+s)^n(s^2)^{2k-n}s^{l-k}.
\endalign$$
Replacing $t$ by $st$ we obtain that
$$\align&[t^n]\l((s+st)^2+s(1-s-st)^2\r)^l
\\=&s^n\sum_{k=0}^l\bi lk\bi{2k}ns^{3k+l-2n}(1+s)^{n-l}
\\=&\sum_{m=0}^{\infty}s^m\sum_{k=0}^l\bi lk\bi{2k}n\bi{n-l}{m+n-3k-l}.
\endalign$$

In view of the above, we immediately get (1.1) for any $m\in\N$
by equating coefficients of $s^m$. \qed

\heading{3. Proof of Theorem 1.2}\endheading

\proclaim{Lemma 3.1} Let $p$ be any prime, and $k,r\in\{0,\ldots,p-1\}$. Then
$$\bi{p-1}k\eq(-1)^k\ (\mo\ p)\
\t{and}\ \bi{p+k+r}{p+r}\eq\bi{p+k+r}r\eq\bi{k+r}r\ (\mo\ p).$$
\endproclaim
\Proof. Clearly
$$\bi{p-1}k=\prod_{0<j\ls k}\f{p-j}j=\prod_{0<j\ls k}\l(\f pj-1\r)\eq(-1)^k\ (\mo\ p)$$
and
$$\align \bi{p+k+r}{p+r}=&\bi{p+k+r}k=\prod_{0<j\ls k}\l(\f pj+1\r)
\times\prod_{0<s\ls r}\f{p+k+s}{p+s}
\\\eq&\prod_{0<s\ls r}\f{p+k+s}s=\bi{p+k+r}r\eq\prod_{0<s\ls r}\f{k+s}s=\bi{k+r}r\ (\mo\ p).
\endalign$$
So we have the desired congruences. \qed

\medskip

\noindent{\it Proof of Theorem 1.2}. In the case $d=p$, (1.4)--(1.6) hold trivially.
Below we assume $d<p$.

(i) Let $m=2(p-1)+d$. Applying (1.2) and (1.3) with $l=p-1$ we obtain that
$$\sum_{k=0}^{p-1}(-1)^{d-k}\bi{p-1}k\bi{m-k}p\bi{2k}{k+d}
=(1-[3\mid m-1])\bi{p-1}{\lceil m/3\rceil}\bi{2\lceil m/3\rceil}p$$
and
$$\sum_{k=0}^{p-1}(-1)^{d-k}\bi{p-1}k\bi{m-k}{p+1}\bi{2k}{k+d}
=(1+[3\mid{m+1}])\bi{p-1}{\lceil m/3\rceil}\bi{2\lceil m/3\rceil}{p+1}.$$

If $d\ls k\ls p-1$, then $0\ls (m-k)-p=p-2-(k-d)<p$ unless $d=0$ and $k=p-1$,
in which case $\bi{2k}{k+d}=\bi{2p-2}{p-1}\eq0\ (\mo\ p)$;
also, $m-k=p\eq0\ (\mo\ p)$ when $m-k-1<p$. Thus,
by applying Lemma 3.1, we have
$$\sum_{k=0}^{p-1}(-1)^{k}\bi{p-1}k\bi{m-k}p\bi{2k}{k+d}
\eq\sum_{k=0}^{p-1}\bi{2k}{k+d}\ (\mo\ p)$$
and
$$\sum_{k=0}^{p-1}(-1)^{k}\bi{p-1}k\bi{m-k}{p+1}\bi{2k}{k+d}
\eq\sum_{k=0}^{p-1}(m-k)\bi{2k}{k+d}\ (\mo\ p).$$

 Let $\ve=(\f{p-d}3)$. Then $p-\ve\eq d\ (\mo\ 3)$. Clearly
$$0<\f m3\ls\l\lceil\f m3\r\rceil=\f{2(p-\ve)+d}3-[3\mid{p-d+1}]\ls\f{m+2}3=\f{2p+d}3<p.$$
If $p\gs 5$ then $\lceil m/3\rceil\gs m/3>p/2$; if $p=3$ then
$\lceil m/3\rceil\gs\lceil 4/3\rceil=2>p/2$.
So
$$0<2\l\lceil\f {m}3\r\rceil-p<2p-p=p$$
unless $p=2$ in which case $\lceil m/3\rceil=1$.
Therefore, with the help of Lemma 3.1,
$$\align&(1-[3\mid{m-1}])\bi{p-1}{\lceil m/3\rceil}\bi{2\lceil m/3\rceil}p
\\\eq&(1-[3\mid{p-d}])(-1)^{\lceil m/3\rceil}
=(1-[3\mid{p-d}])(-1)^{d-[3\mid{p-d+1}]}=(-1)^{d}\ve\pmod{p}
\endalign$$
and
$$\align&(1+[3\mid {m+1}])\bi{p-1}{\lceil m/3\rceil}\bi{2\lceil m/3\rceil}{p+1}
\\\eq&(1+[3\mid{p-d+1}])(-1)^{\lceil m/3\rceil}
\bi{2\lceil m/3\rceil}1[p\not=2]
\\\eq&(1+[3\mid{p-d+1}])(-1)^{d-[3\mid{p-d+1}]}2
\l(\f{2(p-\ve)+d}3-[3\mid{p-d+1}]\r)[p\not=2]
\\\eq&\cases(-1)^d\f 43(1-d)+(-1)^d[p=3]\pmod{p}&\t{if}\ \ve=-1\ (\t{i.e.},\ 3\mid p-d+1),
\\(-1)^d\f23(d-2\ve)+(-1)^d[p=3]\pmod{p}&\t{otherwise}.\endcases
\endalign$$

In view of the above,
$$\sum_{k=0}^{p-1}\bi{2k}{k+d}
\eq(-1)^d(1-[3\mid{m-1}])\bi{p-1}{\lceil m/3\rceil}\bi{2\lceil m/3\rceil}p\eq\ve\pmod{p},$$
and
$$\align &\sum_{k=0}^{p-1}k\bi{2k}{k+d}
\\=&m\sum_{k=0}^{p-1}\bi{2k}{k+d}-\sum_{k=0}^{p-1}(m-k)\bi{2k}{k+d}
\\\eq&m\ve-(-1)^d(1+[3\mid{m+1}])\bi{p-1}{\lceil m/3\rceil}\bi{2\lceil m/3\rceil}{p+1}
\\\eq&\cases(d-2)\ve-\f 43(1-d)-[p=3]=\f{d+2}3-[p=3]&\t{if}\ \ve=-1,
\\(d-2)\ve-\f 23(d-2\ve)-[p=3]=\ve d-\f 23(d+\ve)-[p=3]&\t{otherwise},\endcases
\\\eq&\l([3\mid{p-d}]-\f13\r)(2\ve-d)-[p=3]\ \ (\mo\ p).
\endalign$$
This proves (1.4) and (1.5).

(ii) Our strategy to deduce (1.6) is to compute $S\ \mo\ p^2$ in two different ways, where
$$S=\sum_{k=0}^p(-1)^{k}\bi{p}{k}\bi{2p+d-k}{p}\bi{2k}{k+d}.$$

Observe that $\bi{2p}p\eq 2\ (\mo\ p^2)$. In the case $p\not=2$, this is because
$$\align&\f12\bi{2p}{p}=\bi{2p-1}{p-1}
=\prod_{k=1}^{p-1}\f{2p-k}k=\prod_{k=1}^{p-1}\l(1-\f{2p}k\r)
\\\eq&1-2p\sum_{k=1}^{p-1}\f1k=1-p\sum_{k=1}^{p-1}\l(\f1k+\f1{p-k}\r)\eq1\pmod{p^2}.
\endalign$$
(Moreover, by Wolstenholme's theorem, $\bi{2p-1}{p-1}\eq1\ (\mo\ p^3)$ if $p>3$.)
Therefore
$$\align
S=&\sum_{k=d}^p(-1)^k\bi{p}{k}\bi{2p+d-k}{p}\bi{2k}{k+d}
\\\eq&2(-1)^d\bi{p}{d}-\bi{p+d}{p}\bi{2p}{p+d}
+\sum_{d<k<p}(-1)^k\bi{p}{k}\bi{2k}{k+d}\pmod{p^2}.
\endalign$$
(Note that if $d<k<p$ then $p\mid\bi pk$ and $\bi{p+(p-(k-d))}p\eq1\pmod{p}$.)
If $d\not=0$, then
$$\bi pd=\f pd\bi{p-1}{d-1}\eq(-1)^{d-1}\f pd\ \pmod{p^2}$$
and
$$\align\bi{p+d}p\bi{2p}{p+d}=&\bi{p+d}p\f{2p}{p-d}\bi{2p-1}{p-d-1}
\\=&\f{2p}{p-d}\bi{p+d}p\prod_{0<j<p-d}\l(\f{2p}j-1\r)
\\\eq&\f{2p}{p-d}\bi d0(-1)^{p-d-1}\eq(-1)^{d-1}\f{2p}d\pmod{p^2}.
\endalign$$
Thus
$$\align &S-\sum_{k=d}^{p-1}(-1)^k\bi pk\bi{2k}{k+d}
\\\eq&(-1)^d\bi pd-\bi{p+d}d\bi{2p}{p+d}
\\\eq&\cases -\f pd-(-1)^{d-1}\f{2p}d=(2(-1)^d-1)\f pd\pmod{p^2}&\t{if}\ d\not=0,
\\1-2=-1\pmod{p^2}&\t{if}\ d=0.\endcases
\endalign$$

Set $m=2p+d$. Applying (1.0) with $l=p$ we get that
$$(-1)^dS=\cases\bi{p}{m/3}\bi{2m/3}{p}
&\t{if}\ 3\mid m\ (\t{i.e.,}\ p\eq d\pmod{3}),\\0&\t{otherwise}.\endcases$$
In the case $3\mid m$,
$$\align\bi p{m/3}\bi{2m/3}p=&\f p{m/3}\bi{p-1}{m/3-1}\bi{p+(p+2d)/3}p
\\\eq&\f p{m/3}(-1)^{m/3-1}\eq(-1)^{d-1}\f{3p}m\pmod{p^2}.
\endalign$$
Therefore
$$S\eq\cases-3p/m\pmod{p^2}&\t{if}\ p\eq d\pmod{3},\\0\pmod{p^2}&\t{otherwise}.\endcases$$

Comparing the two congruences for $S\ \mo\ p^2$, we finally obtain that
$$\sum_{k=0}^{p-1}(-1)^k\bi pk\bi{2k}{k+d}
\eq\cases-[3\mid{p-d}]\f{3p}d-(2(-1)^d-1)\f pd\pmod{p^2}&\t{if}\ d\not=0,
\\[p=3]-(-1)\pmod{p^2}&\t{if}\ d=0.\endcases$$
This is equivalent to (1.6) since $\bi pk\eq \f pk(-1)^{k-1}\pmod{p^2}$
for $k=1,\ldots,p-1$.

The proof of Theorem 1.2 is now complete. \qed

\heading{4. Proof of Theorem 1.3}\endheading

\proclaim{Lemma 4.1} Let $r$ be a positive integer, and let $p\gs 4r+7$ be a prime.
Then $\sum_{k=0}^{p-1}\bi{k+r}rC_k$ is congruent to
$$\sum_{k=0}^{\lfloor(r+1-\ve_r)/3\rfloor}
(-1)^{k+\ve_r}\bi{r+2}{3k+1+\ve_r}\bi{k+(r-2+\ve_r)/3}r$$
modulo $p$ with $\ve_r=(\f{p-r-1}3)$.
\endproclaim

\Proof. Let $l=p-r-1$ and $\da\in\{0,1\}$.
Applying (1.1) with $m=2l-\da$ and $n=p$, we obtain that
$$\sum_{k=0}^l(-1)^{k+\da}\bi lk\bi{2l-\da-k}p\bi{2k}{k-\da}
=\sum_{k=0}^l\bi lk\bi{2k}p\bi{r+1}{l+p-\da-3k}.$$
For $k=0,\ldots,l$ it is apparent that
$$\align\bi{l}k=&\bi{p-r-1}k=\f{(p-1)\cdots(p-k-r)}{(k+r)!}
\times\f{(k+1)\cdots(k+r)}{(p-1)\cdots(p-r)}
\\\eq&(-1)^{k+r}(-1)^r\f{(k+1)\cdots(k+r)}{r!}=(-1)^k\bi{k+r}r\ (\mo\ p).
\endalign$$
Thus
$$\align&(-1)^{\da}\sum_{k=0}^l\bi{k+r}r\bi{2l-\da-k}p\bi{2k}{k-\da}
\\\eq&\sum_{k=0}^l(-1)^k\bi{k+r}r\bi{2k}p\bi{r+1}{l+p-\da-3k}\ (\mo\ p).
\endalign$$
If $0\ls k\ls 2l-\da-p=p-2r-2-\da$, then $2l-\da-k\in[p,2p)$ and hence
$$\bi{2l-\da-k}p\eq\bi{2l-\da-k}0=1\ (\mo\ p)$$
by Lemma 3.1.
For $k\in\{0,\ldots,l\}$, clearly
$$r+1\gs l+p-\da-3k\iff 3k\gs 2l-\da\iff k\gs\f{2l-\da}3.$$
If $l\gs k\gs\lceil (2l-\da)/3\rceil$, then
$$2p>2k\gs\f 23(2p-2r-2-\da)=p+\f{p-4r-4-2\da}3>p+\f{p-4r-7}3\gs p$$
and hence $\bi{2k}p\eq\bi{2k}0=1\ (\mo\ p)$ by Lemma 3.1.
Therefore
$$\align&(-1)^\da\sum_{k=0}^{p-2r-2-\da}\bi{k+r}r\bi{2k}{k-\da}
\\\eq&\sum_{\lceil\f{2l-\da}3\rceil\ls k\ls\lfloor\f{l+p}3\rfloor}
(-1)^k\bi{k+r}r\bi{r+1}{l+p-\da-3k}\ (\mo\ p).
\endalign$$

When $p-2r-2\ls k\ls p-2$, we have $2k\gs 2(p-2r-2)>p>k+1$ and hence
$$C_k=\f{(2k)!}{k!(k+1)!}\eq0\ (\mo\ p)\ \ \t{and}\ \ \bi{2k}k=(k+1)C_k\eq0\ (\mo\ p).$$
Note also that $\bi{p-1+r}r=p\cdots(p+r-1)/r!\eq0\ (\mo\ p)$.
So, by the above,
$$\align\sum_{k=0}^{p-1}\bi{k+r}rC_k\eq&\sum_{k=0}^{p-2r-3}\bi{k+r}rC_k
\\\eq&\sum_{k=0}^{p-2r-2}\bi{k+r}r\bi{2k}k-\sum_{k=0}^{p-2r-3}\bi{k+r}r\bi{2k}{k-1}
\\\eq&\sum_{2l\ls 3k\ls l+p}(-1)^k\bi{k+r}r\bi{r+1}{l+p-3k}
\\&+\sum_{2l-1\ls 3k\ls l+p}(-1)^k\bi{k+r}r\bi{r+1}{l+p-1-3k}\ (\mo\ p)
\endalign$$
and hence
$$\sum_{k=0}^{p-1}\bi{k+r}rC_k\eq S(r)\ (\mo\ p),$$
where
$$\align S(r)=&\sum_{2l-1\ls 3k\ls l+p}(-1)^k\bi{k+r}r\bi{r+2}{l+p-3k}
\\=&\sum_{3k\gs 2l-1}(-1)^k\bi{k+r}r\bi{r+2}{3k-2l+1}.
\endalign$$
Clearly
$$\l\lceil\f{2l-1}3\r\rceil=\f{2l-(\f{2l}3)}3=\f{2l+\ve_r}3$$
and thus
$$\align S(r)=&\sum_{k\in\N}(-1)^{k+(2l+\ve_r)/3}\bi{k+(2l+\ve_r)/3+r}r
\bi{r+2}{3k+2l+\ve_r-2l+1}
\\=&\sum_{k\in\N}(-1)^{k+\ve_r}\bi{r+2}{3k+1+\ve_r}\bi{k+(2p+r-2+\ve_r)/3}r
\\\eq&\sum_{k=0}^{\lfloor(r+1-\ve_r)/3\rfloor}
(-1)^{k+\ve_r}\bi{r+2}{3k+1+\ve_r}\bi{k+(r-2+\ve_r)/3}r\ (\mo\ p).
\endalign$$
(Note that $\bi xr$ is a polynomial in $x$ with $p$-adic integer coefficients.)
So we have the desired result. \qed
\medskip

\noindent{\it Proof of Theorem 1.3}.
With the help of (1.14),
$$\align\sum_{k=0}^{p-1}\bi{k+r}rC_{k+d}
=&\sum_{k=d}^{p-1}\bi{k-d+r}rC_k+\sum_{0\ls k<d}\bi{p+k-d+r}rC_{p+k}
\\\eq&\sum_{k=d}^{p-1}\bi{k-d+r}rC_k+2\sum_{0\ls k<d}\bi{k-d+r}rC_{k}\ (\mo\ p).
\endalign$$
By the transformation $(-1)^r\bi{-x}r=\bi{x+r-1}r$ and the Chu-Vandermonde identity,
for any $k\in\{0,\ldots,p-1\}$ we have
$$\align (-1)^r\bi{k-d+r}r=&\bi{d-k-1}r=\sum_{i=0}^r\bi d{r-i}\bi{-k-1}i
\\=&\sum_{i=0}^r\bi d{r-i}(-1)^i\bi{k+i}i.
\endalign$$
Therefore
$$\align&(-1)^r\sum_{k=0}^{p-1}\bi{k+r}rC_{k+d}
\\\eq&\sum_{i=0}^r(-1)^i\bi d{r-i}\sum_{k=0}^{p-1}\bi{k+i}iC_k
+\sum_{0\ls k<d}\bi{d-1-k}rC_k\ \ (\mo\ p).
\endalign$$
So, it suffices to show that $\sum_{k=0}^{p-1}\bi{k+i}iC_k\eq f_i(\ve_i)\ (\mo\ p)$
for all $i=0,\ldots,p-1$. As $r$ is an arbitrarily chosen element of $\{0,\ldots,p-1\}$,
below we only need to show the congruence
$$\sum_{k=0}^{p-1}\bi{k+r}rC_k\eq f_r(\ve_r)\ \ (\mo\ p).\tag4.1$$

To prove (4.1) we further extend the idea in the proofs of (1.4) and (1.5).

Let $\da\in\{0,1\}$.
Applying (1.1) with $l=p-1$, $m=2p-1-\da$ and $n=p+r$ we get that
$$(-1)^{\da+1}S_{\da}=\sum_{k=0}^{p-1}\bi{p-1}k\bi{2k}{p+r}\bi{r+1}{2p-\da+r-3k},$$
where
$$S_{\da}=\sum_{k=0}^{p-1}(-1)^{k}\bi{p-1}k\bi{2p-1-\da-k}{p+r}\bi{2k}{k+1-\da}.$$

By Lemma 3.1, $(-1)^k\bi{p-1}k\eq1\ (\mo\ p)$ for all $k=0,\ldots,p-1$, and
$\bi{K}{p+r}\eq\bi Kr\ (\mo\ p)$ for any integer $K\in[p+r,2p+r)$.
Thus
$$S_{\da}\eq\sum_{0\ls k<p-r-\da}\bi{2p-1-\da-k}r\bi{2k}{k+1-\da}\ \ (\mo\ p)$$
and
$$\align&\sum_{k=0}^{p-1}\bi{p-1}k\bi{2k}{p+r}\bi{r+1}{2p-\da+r-3k}
\\\eq&\sum_{\lceil\f{p+r}2\rceil\ls k<p}(-1)^k\bi{2k}r\bi{r+1}{2p-\da+r-3k}\ \ (\mo\ p).
\endalign$$
Therefore
$$\aligned&(-1)^{\da+1}\sum_{0\ls k<p-r-\da}\bi{2p-1-\da-k}r\bi{2k}{k+1-\da}
\\\eq&\sum_{\lceil\f{p+r}2\rceil\ls k<p}(-1)^k\bi{2k}r\bi{r+1}{2p-\da+r-3k}\ (\mo\ p).
\endaligned\tag4.2$$

If $r\in\{1,\ldots,p-1\}$ then
$$\bi{2p-2-(p-1-r)}r=\bi{p-1+r}r=\f{(p+r-1)!}{(p-1)!r!}\eq0\ \ (\mo\ p);$$
if $r=0$ then
$$\bi{2(p-1-r)}{p-1-r}=\f{(p+p-2)!}{(p-1)!(p-1)!}\eq0\ \ (\mo\ p).$$
Thus, when $k=p-1-r$ we have
$$\bi{2p-2-k}r\bi{2k}k\eq0\ \ (\mo\ p).$$
In view of this and (4.2),
$$\align &\sum_{k=0}^{p-1-r}\bi{2p-1-k}rC_k
\\=&\sum_{k=0}^{p-1-r}\bi{2p-1-k}r\(\bi{2k}k-\bi{2k}{k+1}\)
\\=&-\sum_{k=0}^{p-1-r}\bi{2p-1-k}r\bi{2k}{k+1}
\\&+\sum_{k=0}^{p-1-r}\(\bi{2p-2-k}r+\bi{2p-2-k}{r-1}\)\bi{2k}k
\\\eq&\sum_{\lceil\f{p+r}2\rceil\ls k<p}(-1)^k\bi{2k}r\(\bi{r+1}{2p+r-3k}+\bi{r+1}{2p-1+r-3k}\)
\\&+\sum_{\lceil\f{p+r-1}2\rceil\ls k<p}(-1)^k\bi{2k}{r-1}\bi r{2p-2+r-3k}\ (\mo\ p).
\endalign$$
Since
$$\bi{2p-1-k}r=\prod_{0<s\ls r}\f{2p-k-s}s
\eq(-1)^r\prod_{0<s\ls r}\f{k+s}s=(-1)^r\bi{k+r}r\ (\mo\ p)$$
for every $k=0,\ldots,p-1-r$, we have
$$\align&(-1)^r\sum_{k=0}^{p-1-r}\bi{k+r}rC_k
\\\eq&\sum_{\lceil\f{p+r}2\rceil\ls k<p}(-1)^k\bi{2k}r\bi{r+2}{2p+r-3k}
\\&+\sum_{\lceil\f{p+r-1}2\rceil\ls k<p}(-1)^k\bi{2k}{r-1}\bi r{2p-2+r-3k}\ (\mo\ p).
\endalign$$
When $0\ls 2p+r-3k\ls r+2$ (i.e., $2p-2\ls 3k\ls 2p+r$), if $k<(p+r)/2$ then
$p-1<(4p-4)/3\ls 2k\ls p+r-1$ and hence
$$\bi{2k}r=\f{2k(2k-1)\cdots(2k-r+1)}{r!}\eq0\ (\mo\ p).$$
Similarly, when $0\ls 2p-2+r-3k\ls r$
(i.e., $2p-2\ls 3k\ls 2p-2+r$), if $k<(p+r-1)/2$ then $p\ls 2k\ls p+r-2$
and hence
$$\bi{2k}{r-1}=\f{2k(2k-1)\cdots(2k-r+2)}{(r-1)!}\eq0\ (\mo\ p).$$
Therefore
$$\align&(-1)^r\sum_{k=0}^{p-1-r}\bi{k+r}rC_k
\\\eq&\sum_{2p-2\ls 3k\ls 2p+r}(-1)^k\bi{2k}r\bi{r+2}{2p+r-3k}
\\&+\sum_{2p-2\ls 3k\ls 2p-2+r}(-1)^k\bi{2k}{r-1}\bi r{2p-2+r-3k}
\\\eq&\sum_{k\gs h}(-1)^k\(\bi{2k}r\bi{r+2}{3k-2p+2}
+\bi{2k}{r-1}\bi r{3k-2p+2}\)
\ (\mo\ p),
\endalign$$
where
$$h=\l\lceil\f{2p-2}3\r\rceil=\f{2p-1+\ve_1}3.$$
(Recall that $\ve_1=(\f{p-2}3)\eq p-2\ (\mo\ 3)$.)

If $p-r\ls k\ls p-1$, then $k+1\ls p\ls k+r$ and hence
$$\bi{k+r}r=\prod_{s=1}^r\f{k+s}s\eq0\ \ (\mo\ p).$$
Note also that $2h\eq2(\ve_1-1)/3+[p=3]\ (\mo\ p)$.
Thus, by the above we have
$$\align&(-1)^r\sum_{k=0}^{p-1}\bi{k+r}rC_k
\\\eq&\sum_{k\in\N}(-1)^{k+h}\(\bi{2k+2h}r\bi{r+2}{3k+3h-2p+2}
+\bi{2k+2h}{r-1}\bi r{3k+3h-2p+2}\)
\\\eq& \Psi_r(p\ \mo\ 3)\ \ (\mo\ p),
\endalign$$
where
$$\align \Psi_r(p\ \mo\ 3)=&
\sum_{k\in\N}(-1)^{k+\ve_1-1}\bi{2k+2(\ve_1-1)/3+[p=3]}r\bi{r+2}{3k+1+\ve_1}
\\&+\sum_{k\in\N}(-1)^{k+\ve_1-1}\bi{2k+2(\ve_1-1)/3+[p=3]}{r-1}\bi r{3k+1+\ve_1}.
\endalign$$
(Note that both $\ve_1$ and $[p=3]$ only depend on $p$ mod 3.)

As
$$\Psi_0(p\ \mo\ 3)=(-1)^{\ve_1-1}\bi2{1+\ve_1}=f_0(\ve_0),$$
(4.1) holds when $r=0$. If $p=3$ then
$$-\Psi_1(p\ \mo\ 3)=-3\eq 0=f_1(\ve_1)\ (\mo\ 3)
\ \ \t{and}\ \ \Psi_2(p\ \mo\ 3)=1=f_2(\ve_2).$$
So (4.1) is also valid in the case $p=3$.

Below we assume that $r\not=0$ and $p\not=3$.
Recall that
$$\sum_{k=0}^{p-1}\bi{k+r}rC_k\eq(-1)^r\Psi_r(p\ \mo\ 3)\ \ (\mo\ p).$$
If $p'\gs 4r+7$ is a prime with $p'\eq p\ (\mo\ 3)$, then
$$(-1)^r\Psi_r(p\ \mo\ 3)=(-1)^r\Psi_r(p'\ \mo\ 3)
\eq\sum_{k=0}^{p'-1}\bi{k+r}rC_k\eq f_r(\ve_r)\ (\mo\ p')$$
with the help of Lemma 4.1.
By Dirichlet's theorem (cf. [IR, p.\,251]),
there are infinitely many primes $p'$ with $p'\eq p\ (\mo\ 3)$.
So we must have $(-1)^r\Psi_r(p\ \mo\ 3)=f_r(\ve_r)$ and hence (4.1) follows.
We are done. \qed
\medskip

\heading{5. Proof of Theorem 1.4}\endheading

In this section we let $p$ be an arbitrary prime greater than $d$.

In view of (1.13),
$$\sum_{k=0}^{p-1}C_{k+d}=C_d\sum_{k=0}^{p-1}\bi{2k}k
+\sum_{j=1}^{d+1}C_{d,j}\sum_{k=0}^{p-1}\bi{2k}{k+j}$$
and
$$\sum_{k=0}^{p-1}kC_{k+d}=C_d\sum_{k=0}^{p-1}k\bi{2k}k
+\sum_{j=1}^{d+1}C_{d,j}\sum_{k=0}^{p-1}k\bi{2k}{k+j}.$$
Combining these with Theorem 1.2, we immediately get the congruences
$$\sum_{k=0}^{p-1}C_{k+d}
\eq\l(\f p3\r)C_d+\sum_{j=1}^{d+1}\l(\f{p-j}3\r)C_{d,j}\ (\mo\ p)\tag5.1$$
and
$$\sum_{k=0}^{p-1}kC_{k+d}\eq-\l(\f p3\r)\f23C_d
+\sum_{j=1}^{d+1}C_{d,j}
\l([3\mid{p-j}]-\f13\r)\l(2\l(\f{p-j}3\r)-j\r)\ (\mo\ p).
\tag5.2$$
(It is easy to check that $C_d+\sum_{j=1}^{d+1}C_{d,j}=0$
if $p=3$ (and hence $d\in\{0,1,2\}$).)

By (1.16) and (5.1),
$$\f{3(\f p3)-1}2+\sum_{0\ls k<d}C_k\eq
\l(\f p3\r)C_d+\sum_{j=1}^{d+1}\l(\f{p-j}3\r)C_{d,j}\ (\mo\ p).$$
If $p$ is congruent to $1$ or $2$ modulo $3$, this gives
$$1+\sum_{0\ls k<d}C_k\eq
C_d-\sum_{j=1}^{d+1}\l(\f{j-1}3\r)C_{d,j}\ (\mo\ p)$$
and
$$-2+\sum_{0\ls k<d}C_k\eq
-C_d-\sum_{j=1}^{d+1}\l(\f{j-2}3\r)C_{d,j}\ (\mo\ p)$$
respectively. Note that both sides of these two congruences
are independent of $p$. Thus we have the first equality in (1.19)
since the residue classes $1\,(\mo\ 3)$
and $2\,(\mo\ 3)$ both contain infinitely many primes by Dirichlet's theorem.
The second equality in (1.19) also holds because
$$\align&\sum_{\da=0}^1\((1-2\da)\sum_{0\ls k<d}C_k
+(-1)^{\da}\sum_{i=0}^d\l(\f{i-\da}3\r)C_{d,i+1}+1+\da\)
\\=&\sum_{i=0}^d\l(\l(\f i3\r)-\l(\f{i-1}3\r)\r)C_{d,i+1}+1+2
=\sum_{i=0}^d(1-3[3\mid i+1])C_{d,i+1}+3.
\endalign$$

Observe that
$$\align&\sum_{k=0}^{p-1}kC_{k+d}+(d+1)\sum_{k=0}^{p-1}C_{k+d}
\\=&\sum_{k=0}^{p-1}\bi{2k+2d}{k+d}=\sum_{j=-d}^d\bi{2d}{d-j}\sum_{k=0}^{p-1}\bi{2k}{k+j}
\\=&\bi{2d}d\sum_{k=0}^{p-1}\bi{2k}k+2\sum_{0<j\ls d}\bi{2d}{d-j}\sum_{k=0}^{p-1}\bi{2k}{k+j}.
\endalign$$
Applying Theorem 1.2 and (5.1), we get that
$$\align\sum_{k=0}^{p-1}kC_{k+d}
\eq&\l(\f p3\r)\bi{2d}d+2\sum_{0<j\ls d}\l(\f{p-j}3\r)\bi{2d}{d-j}
\\&-(d+1)\l(\f p3\r)C_d-(d+1)\sum_{j=1}^{d+1}\l(\f{p-j}3\r)C_{d,j}
\\\eq&\sum_{j=1}^{d+1}\l(\f{p-j}3\r)\(2\bi{2d}{d-j}-(d+1)C_{d,j}\)
\ (\mo\ p). \endalign$$
Comparing this with (1.17) we obtain the identity (1.20) by applying Dirichlet's theorem.

It follows from (5.1) and (5.2) that
$$\sum_{k=0}^{p-1}\l(k+\f 23\r)C_{k+d}
\eq\sum_{j=1}^{d+1}\l(\f13-[3\mid{p-j}]\r)jC_{d,j}\ (\mo\ p).\tag5.3$$
On the other hand, by (1.16) and (1.17) we have
$$\align\sum_{k=1}^{p-1}\l(k+\f23\r)C_{k+d}
\eq &\f{d+1}2\l(1-\l(\f p3\r)\r)-\l(\f p3\r)d
-\sum_{k=0}^dkC_{d-k}
\\&+\f 23\cdot\f{3(\f p3)-1}2+\f 23\sum_{0\ls k<d}C_k\ \ (\mo\ p).
\endalign$$
Comparing this with (5.3) we finally get (1.21) by applying Dirichlet's theorem.

The proof of Theorem 1.4 is now complete.

\Ack. The authors are grateful to Prof. C. Krattenthaler, H. Wilf
and the referees for their helpful comments.


\widestnumber\key{PWZ}

 \Refs

\ref\key AAR\by G. E. Andrews, R. Askey and R. Roy\book Special Functions
\publ Cambridge Univ. Press, Cambridge, 1999\endref

\ref\key B\by W. Bailey\paper Products of generalized hypergeometric series
\jour Proc. London Math. Soc.\vol 28\yr 1928\pages 242--254\endref

\ref\key C\by D. Callan\paper A combinatorial proof of Sun's
``curious" identity \jour Integers\vol 4\yr 2004\pages A5, 6 pp. (electronic)\endref

\ref\key GJ\by I. P. Goulden and D. M. Jackson
\book Combinatorial Enumeration\publ John Wiley $\&$ Sons, New York, 1983\endref

\ref\key GKP\by R. L. Graham, D. E. Knuth and O. Patashnik
 \book Concrete Mathematics\publ 2nd ed., Addison-Wesley, New York\yr 1994\endref

\ref\key IR\by K. Ireland and M. Rosen
\book A Classical Introduction to Modern Number Theory
{\rm (Graduate texts in math.; 84), 2nd ed.}
\publ Springer, New York, 1990\endref

\ref\key KR\by C. Krattenthaler and K. S. Rao
\paper Automatic generation of hypergeometric identities by the beta integral method
\jour J. Comput. Appl. Math.\vol 160\yr 2003\pages 159--173\endref

\ref\key PWZ\by M. Petkov\v sek, H. S. Wilf and D. Zeilberger
\book $A=B$\publ A K Peters, Wellesley, 1996\endref

\endRefs
\enddocument